\documentstyle [12pt,epsf, graphicx, amssymb]{article}

\begin{document}
\def\l{\lambda}
\def\m{\mu}
\def\a{\alpha}
\def\b{\beta}
\def\g{\gamma}
\def\d{\delta}
\def\e{\epsilon}
\def\o{\omega}
\def\O{\Omega}
\def\v{\varphi}
\def\t{\theta}
\def\r{\rho}
\def\bs{$\blacksquare$}
\def\bp{\begin{proposition}}
\def\ep{\end{proposition}}
\def\bt{\begin{th}}
\def\et{\end{th}}
\def\be{\begin{equation}}
\def\ee{\end{equation}}
\def\bl{\begin{lemma}}
\def\el{\end{lemma}}
\def\bc{\begin{corollary}}
\def\ec{\end{corollary}}
\def\pr{\noindent{\bf Proof: }}
\def\note{\noindent{\bf Note. }}
\def\bd{\begin{definition}}
\def\ed{\end{definition}}
\def\C{{\mathbb C}}
\def\P{{\mathbb P}}
\def\Z{{\mathbb Z}}
\def\d{{\rm d}}
\def\deg{{\rm deg\,}}
\def\deg{{\rm deg\,}}
\def\arg{{\rm arg\,}}
\def\min{{\rm min\,}}
\def\max{{\rm max\,}}

\newcommand{\norm}[1]{\left\Vert#1\right\Vert}
\newcommand{\abs}[1]{\left\vert#1\right\vert}

\newcommand{\set}[1]{\left\{#1\right\}}
\newcommand{\setb}[2]{ \left\{#1 \ \Big| \ #2 \right\} }

\newcommand{\IP}[1]{\left<#1\right>}
\newcommand{\Bracket}[1]{\left[#1\right]}
\newcommand{\Soger}[1]{\left(#1\right)}

\newcommand{\Integer}{\mathbb{Z}}
\newcommand{\Rational}{\mathbb{Q}}
\newcommand{\Real}{\mathbb{R}}
\newcommand{\Complex}{\mathbb{C}}

\newcommand{\eps}{\varepsilon}
\newcommand{\To}{\longrightarrow}
\newcommand{\varchi}{\raisebox{2pt}{$\chi$}}

\newcommand{\E}{\mathbf{E}}
\newcommand{\Var}{\mathrm{var}}

\def\squareforqed{\hbox{\rlap{$\sqcap$}$\sqcup$}}
\def\qed{\ifmmode\squareforqed\else{\unskip\nobreak\hfil
\penalty50\hskip1em\null\nobreak\hfil\squareforqed
\parfillskip=0pt\finalhyphendemerits=0\endgraf}\fi}

\renewcommand{\th}{^{\mathrm{th}}}
\newcommand{\Dif}{\mathrm{D_{if}}}
\newcommand{\Difp}{\mathrm{D^p_{if}}}
\newcommand{\GHF}{\mathrm{G_{HF}}}
\newcommand{\GHFP}{\mathrm{G^p_{HF}}}
\newcommand{\f}{\mathrm{f}}
\newcommand{\fgh}{\mathrm{f_{gh}}}
\newcommand{\T}{\mathrm{T}}
\newcommand{\K}{^\mathrm{K}}
\newcommand{\PghK}{\mathrm{P^K_{f_{gh}}}}
\newcommand{\Dig}{\mathrm{D_{ig}}}
\newcommand{\for}{\mathrm{for}}
\newcommand{\End}{\mathrm{end}}

\newtheorem{th}{Theorem}[section]
\newtheorem{lemma}{Lemma}[section]
\newtheorem{definition}{Definition}[section]
\newtheorem{corollary}{Corollary}
\newtheorem{proposition}{Proposition}[section]

\begin{titlepage}

\begin{center}

\topskip 5mm

{\LARGE{\bf {Remez-Type Inequality for Smooth Functions}}} \vskip 8mm

{\large {\bf Y. Yomdin}}

\vspace{6 mm}
\end{center}

{Department of Mathematics, The Weizmann Institute of Science,
Rehovot 76100, Israel}

\vspace{6 mm}

{e-mail: yosef.yomdin@weizmann.ac.il}

\vspace{1 mm}

\vspace{6 mm}
\begin{center}

{ \bf Abstract}
\end{center}

{\small The classical Remez inequality bounds the maximum of the absolute value of a polynomial $P(x)$ of degree $d$ on $[-1,1]$
through the maximum of its absolute value on any subset $Z$ of positive measure in $[-1,1]$. Similarly, in several variables the
maximum of the absolute value of a polynomial $P(x)$ of degree $d$ on the unit ball $B^n \subset {\mathbb R}^n$ can be bounded
through the maximum of its absolute value on any subset $Z\subset Q^n_1$ of positive $n$-measure $m_n(Z)$. In \cite{Yom} a stronger
version of Remez inequality was obtained: the Lebesgue $n$-measure $m_n$ was replaced by a certain geometric quantity $\omega_{n,d}(Z)$
satisfying $\omega_{n,d}(Z)\geq m_n(Z)$ for any measurable $Z$. The quantity $\omega_{n,d}(Z)$ can be effectively estimated in terms of 
the metric entropy of $Z$ and it may be nonzero for discrete and even finite sets $Z$.

In the present paper we extend Remez inequality to functions of finite smoothness. This is done by combining the result of \cite{Yom}
with the Taylor polynomial approximation of smooth functions. As a consequence we obtain explicit lower bounds in some examples in the 
Whitney problem of a $C^k$-smooth extrapolation from a given set $Z$, in terms of the geometry of $Z$.}

\vspace{1 mm}
\begin{center}
------------------------------------------------
\vspace{1 mm}
\end{center}
This research was supported by the ISF, Grant No. 639/09, and by the Minerva Foundation.

\end{titlepage}

\newpage

\section{Introduction}
\setcounter{equation}{0}

The classical Remez inequality (\cite{Rem}, see also \cite{Erd}) reads as follows:

\medskip

\bt\label{Remez.1}
Let $P(x)$ be a polynomial of degree $d$. Then for any measurable $Z\subset [-1,1]$
\be \label{Remez.ineq.1}
\max_{[-1,1]} \vert P(x) \vert \leq T_d({{4-m}\over {m}})\max_Z \vert P(x) \vert ,
\ee
where $m=m_1(Z)$ is the Lebesgue measure of $Z$ and $T_d(x)=cos(d \ arccos(x))$ is the $d$-th Chebyshev polynomial.
\et

In several variables a generalization of Theorem \ref{Remez.1} was obtained in \cite{Bru.Gan}:

\bt\label{Remez.n}
Let ${\cal B}\subset {\mathbb R}^n$ be a convex body and let $\O\subset {\cal B}$ be a measurable set. Then for any real
polynomial $P(x)=P(x_1,\dots,x_n)$ of degree $d$ we have
\be\label{Remez.ineq.n}
\sup_{{\cal B}}\vert P \vert \leq T_d ({{1+(1-\lambda)^{1\over n}}\over {1-(1-\lambda)^{1\over n}}}) \sup_{\O}\vert P \vert.
\ee
Here $\lambda= {{m_n(\O)}\over {m_n({\cal B})}},$ with $m_n$ being the Lebesgue measure on ${\mathbb R}^n$. This inequality
is sharp and for $n=1$ it coincides with the classical Remez inequality.
\et

It is clear that Remez inequality of Theorems \ref{Remez.1} and \ref{Remez.n} cannot be verbally extended to
smooth functions: such function $f$ may be identically zero on any given closed set $Z$, and non-zero elsewhere. In the present paper
we show that adding a ``remainder term'' (expressible through the bounds on the derivatives of $f$) provides a generalization
of the Remez inequality to smooth functions. Our main goal is to study the interplay between the geometry of the ``sampling set''
$Z$, the bounds on the derivatives of $f$, and the bounds for the extension of $f$ from $Z$ to the ball $B^n$ of radius $1$ centered
at the origin in ${\mathbb R}^n$. To state our main ``general'' result we need some definitions:

\bd\label{Remez.constant}
For a set $Z\subset B^n \subset {\mathbb R}^n$ and for each $d\in {\mathbb N}$ the Remez constant $R_d(Z)$ is the minimal $K$
for which the inequality $\sup_{B^n}\vert P \vert \leq K \sup_{Z}\vert P \vert$ is valid for any real polynomial
$P(x)=P(x_1,\dots,x_n)$ of degree $d$.
\ed
For some $Z$ the Remez constant $R_d(Z)$ may be equal to $\infty$. In fact, $R_d(Z)$ is infinite if and only if $Z$ is contained
in the set of zeroes $Y_P=\{x\in {\mathbb R}^n, \ | \ P(x)=0\}$ of a certain polynomial $P$ of degree $d$. See \cite{Bru.Yom} for
a detailed discussion.

\bd\label{Appr.Error}
Let $f: B^n \rightarrow {\mathbb R}$ be a $k$ times continuously differentiable function on $B^n$. For $d=0,1,\ldots,$
the approximation error $E_d(f)$ is the minimum over all the polynomials $P(x)$ of degree $d$ of the absolute deviation
$M_0(f-P)= \max_{x\in B^n}|f(x)-P(x)|$.
\ed

\bt\label{Smooth.Remez}
Let $f: B^n \rightarrow {\mathbb R}$ be a $k$ times continuously differentiable function on $B^n$, and let a subset
$Z\subset B^n$ be given. Put $L=\max_{x\in Z}|f(x)|$. Then

\be\label{Smooth.Remez.ineq}
\max_{x\in B^n} \vert f(x) \vert \leq \inf_d [R_d(Z)(L+E_d(f))+E_d(f)].
\ee
\et
\pr
Let for a fixed $d$ $P_d(x)$ be the polynomial of degree $d$ for which the best approximation of $f$ is achieved:
$E_d(f) = \max_{x\in B^n}|f(x)-P(x)|$. Then $\max_{x\in Z}|P(x)| \leq L+E_d(f)$. By definition of the Remez constant $R_d(Z)$
we have $\max_{x\in B^n}|P(x)| \leq R_d(Z)(L+E_d(f))$. Returning to $f$ we get
$\max_{x\in B^n} \vert f(x) \vert \leq R_d(Z)(L+E_d(f))+E_d(f).$ Since this is true for any $d$, we finally obtain
$\max_{x\in B^n} \vert f(x) \vert \leq \inf_d [R_d(Z)(L+E_d(f))+E_d(f)].$ $\square$

\medskip

In this paper we produce, based on Theorem \ref{Smooth.Remez}, explicit Remez-type bounds for smooth functions in some typical
situations.

\section{Bounding $R_d(Z)$ via Metric entropy}\label{entr.bounds}
\setcounter{equation}{0}

It is well known that the inequality of the form (\ref{Remez.ineq.1}) or (\ref{Remez.ineq.n})
may be true also for some sets $Z$ of measure zero and even for certain discrete or finite sets $Z$. Let us mention here only
a couple of the most relevant results in this direction: in \cite{Cop.Riv,Rah,Yom.Zah} such inequalities are provided
for $Z$ being a regular grid in $[-1,1]$. In \cite{Fav} discrete sets $Z\subset [-1,1]$ are studied. In this last paper the
invariant $\phi_Z(d)$ is defined and estimated in some examples, which is the best constant in the Remez-type inequality of
degree $d$ for the couple $(Z\subset [-1,1])$.

\smallskip

In \cite{Yom} (see also \cite{Bru}) a strengthening of Remez inequality was obtained: the Lebesgue $n$-measure $m_n$ was replaced
by a certain geometric quantity $\omega_{n,d}(Z),$ defined in terms of the metric entropy of $Z$, and satisfying
$\omega_{n,d}(Z)\geq m_n(Z)$ for any measurable $Z\subset Q^n_1$. So we have the following proposition, which combines the result
of Theorem 3.3 of \cite{Yom} with the well-known bound for Chebyshev polynomials (see \cite{Erd}):

\bp\label{R.Const.Entr}
For each $Z\subset B^n$ and for any $d$ the Remez constant $R_{n,d}(Z)$ satisfies

\be\label{R.Const.Entr.eq}
R_{n,d}(Z) \leq T_d ({{1+(1-\lambda)^{1\over n}}\over {1-(1-\lambda)^{1\over n}}})\leq ({4n\over \lambda})^d,
\ee
where $\lambda = \omega_{n,d}(Z)$.
\ep
In what follows we shall omit the dimension $n$ from the notations for $\omega_d(Z)=\omega_{n,d}(Z)$. It was shown in \cite{Yom}
that in many cases (but not always!) the bound of Proposition \ref{R.Const.Entr} is pretty sharp. In the present paper we recall
the definition of $\omega_d(Z)$ and estimate this quantity in several typical cases, stressing the setting where $Z$ is fixed,
while $d$ changes.

\subsection{Definition and properties of $\omega_d(Z)$}

To define $\omega_d(Z)$ let us recall that the covering number $M(\e,A)$ of a metric space $A$ is the minimal number of closed
$\e$-balls covering $A$. Below $A$ will be subsets of ${\mathbb R}^n$ equipped with the $l^\infty$ metric. So the $\e$-balls in
this metric are the cubes $Q^n_\e$.

For a polynomial $P$ on ${\mathbb R}^n$ let us consider the sub-level set $V_\rho(P)$ defined by $V_\rho(P)=\{x\in B^n,
\vert P(x) \vert \leq \rho \}$. The following result is proved in (\cite{Vit1}):

\bt\label{Vit} ({\bf Vitushkin's bound}) For $V=V_\rho(P)$ as above
\be\label{Vit.eq}
M(\e,V)\leq \sum_{i=0}^{n-1}C_i(n,d)({1\over \e})^i+ m_n(V)({1\over{\e}})^n,
\ee
with $C_i(n,d)=C'_i(n)(2d)^{(n-i)}$. For $n=1$ we have $M(\e,V)\leq d + \mu_1(V)({1\over{\e}}),$ and for $n=2$ we have
$$
M(\e,V)\leq (2d-1)^2 + 8d({1\over \e}) + \mu_2(V)({1\over{\e}})^2.
$$
\et
For $\e>0$ we denote by $M_{n,d}(\e)$ (or shortly $M_d(\e)$) the polynomial of degree $n-1$ in ${1\over \e}$ as appears in
(\ref{Vit.eq}):
\be\label{Md}
M_d(\e)=\sum_{i=0}^{n-1}C_i(n,d)({1\over \e})^i.
\ee
In particular,
$$
M_{1,d}(\e)=d, \ M_{2,d}(\e)=(2d-1)^2 + 8d({1\over \e}).
$$
Now for each subset $Z\subset B^n$ (possibly discrete or finite) we introduce the quantity $\omega_d(Z)$ via the following
definition:

\bd\label{omegad}
Let $Z$ be a subset in $B^n\subset {\mathbb R}^n$. Then $\omega_d(Z)$ is defined as
\be\label{omga.eq}
\omega_d(Z) = \sup_{\e>0} {\e}^n[M(\e,Z)- M_d(\e)].
\ee
\ed
The following results are obtained in \cite{Yom}:

\bp\label{prop.o}
The quantity $\o_d(Z)$ for $Z\subset B^n$ has the following properties:

\smallskip

1. For a measurable $Z$ \ \ $\o_d(Z)\geq m_n(Z).$

\smallskip

2. For any set $Z\subset B^n$ the quantities $\o_d(Z)$ form a non-increasing sequence in $d$.

\smallskip

3. For a set $Z$ of Hausdorff dimension $n-1$, if the Hausdorff $n-1$ measure of $Z$ is large enough with respect to $d$, then
$\o_d(Z)$ is positive.

\smallskip

4. Let $G_s=\{x_1=-1,x_2,\dots,x_s=1\}$ be a regular grid in $[-1,1]$. Then $\o_d(G_s)={{2(s-d)}\over {s-1}}$.

\smallskip

Let $Z_r=\{1,{1\over {2^r}},{1\over {3^r}},\dots,{1\over {k^r}},\dots\}$. In this case
$\o_d(Z_r)\asymp {{r^r}\over {(r+1)^{r+1}}} \ {1\over {d^r}}.$

\smallskip

Let $Z(q)=\{1,q,q^2,q^3,\dots,q^m,\dots\}, \ 0< q < 1$. Then $\o_d(Z(q)) \asymp {{q^d}\over {\log ({1\over q})}}$.
\ep
We need the following result, which, although in the direction of the results in \cite{Yom}, was not proved there explicitly.
Let $S$ be a connected smooth curve in $B^2$ of the length $\sigma$. Define $\e_0$ as the maximal $\e$ such that for each
$\delta \leq \e$ we have $M(\delta, S)\geq {{l(S)}\over {2\delta}}$. The parameter $\e_0$ is a kind of ``injectivity radius'' of
the curve $S$, and for any curve of length $\sigma$ inside the unit ball $B^2$ it cannot be larger than ${1\over \sigma}$.
Write $\e_0$ as $\e_0={1\over {l\sigma}}, \ l\geq 1$. The computation below essentially compares the length of $S$ with the maximal
possible length of an algebraic curve of degree $d$ inside $B^2$, which is of order $d$. So it is convenient for any given $d$
to write $\sigma$ as $\sigma=md$.

\bp\label{curve}

In the notations above, $\o_d(S)$ satisfies

\be\label{curve.eq}
\o_d(S)\geq {1\over {2l}}(1-{{24}\over m}).
\ee
In particular, for the length of $S$ larger than $24d$, \ \ $\o_d(S)$ is strictly positive.
\ep
\pr
By definition,
$$
\o_d(S)=\sup_\e \e^2[M(\e,S)-M_d(\e)]=\sup_\e \e^2[M(\e,S)-(2d-1)^2-8d({1\over \e})].
$$
Substituting here $\e_0={1\over {l\sigma}}$ we get

$$
\o_d(S) \geq ({1\over {lmd}})^2[{{l(md)^2}\over 2}-(2d-1)^2-8lmd^2]=
$$
$$
={1\over {2l}}(1-{2\over m}[({{2d-1}\over d})^2+8])\geq {1\over {2l}}(1-{{24}\over m}).
$$
In particular, for $m>24$, i.e. for the length of $S$ larger than $24d$, the quantity $\o_d(S)$ is strictly positive. $\square$

\section{Bounding Smooth Functions}\label{Smooth.f}
\setcounter{equation}{0}

Let $f: B^n \rightarrow {\mathbb R}$ be a $k$ times continuously differentiable function on $B^n$. For $l=0,1,\ldots,k$ put
$M_l(f)=\max_{B^n} \Vert d^l f \Vert,$ where the norm of the $l$-th differential of $f$ is defined as the sum of the absolute
values of all the partial derivatives of $f$ of order $l$. To simplify notations, we shall not make specific assumptions on the
continuity modulus of the last derivative $d^k f$. Now we use Taylor polynomials of an appropriate degree between $0$ and $k-1$ in
order to bound from above the approximation error $E_d(f),\ d=0,1,\ldots,k$. Applying Theorem \ref{Smooth.Remez}, we obtain
the following result:

\bp\label{Taylor1}

Let $f: B^n \rightarrow {\mathbb R}$ be a $k$ times continuously differentiable function on $B^n$, with
$M_l(f)=\max_{B^n} \Vert d^l f \Vert, \ l=0,1,\ldots,k,$ and let a subset $Z\subset B^n$ be given. Put $L=\max_{x\in Z}|f(x)|$.
Then

\be\label{Smooth.Remez.ineq1}
M_0(f)=\max_{x\in B^n} \vert f(x) \vert \leq \min_{d=0,1,\ldots,k-1} [R_d(Z)(L+E^T_d(f))+E^T_d(f)],
\ee
where $E^T_d(f)={1\over {(d+1)!}}M_{d+1}(f)$ is the Taylor remainder term of $f$ of degree $d$ on the unit ball $B^n$.
\ep
\pr
We restrict infinum in Theorem \ref{Smooth.Remez} to a smaller set of $d$'s, and replace $E_d(f)$ with a larger quantity $E^T_d(f)$. 
$\square$.

\smallskip

In general we cannot get an explicit answer for the minimum in Proposition \ref{Taylor1}, unless we add more specific assumptions
on the set $Z$ and the sequence $M_d(f)$. However, this proposition provides an explicit and rather sharp information in the case where
the set $Z$ is ``small''. Let us pose the following question: for a fixed $s=1,\ldots,k-1$ and a given set $Z\subset B^n$ is it possible
to bound $M_0(f)=\max_{x\in B^n} \vert f(x) \vert$ through $L=\max_{x\in Z}|f(x)|$ and $M_{s+1}(f)$ only, without knowing bounds on the
derivatives $d^l(f), \ l\leq s$?

\bp\label{smallZ}
If $R_s(Z)<\infty$ then $M_0(f)\leq R_s(Z)(L+E^T_s(f))+E^T_s(f)$ with $E^T_s(f)={1\over {(s+1)!}}M_{s+1}(f)$. If $R_s(Z)=\infty$ then $M_0(f)$
cannot be bounded in terms of $L$ and $M_l(f), \ l\geq s+1.$
\ep
\pr
In case $R_s(Z)<\infty$ the required bound is obtained by restricting the minimization in (\ref{Smooth.Remez.ineq1}) to $d=s$ only.
If $R_s(Z)=\infty$ then already polynomials of degree $s$ vanishing on $Z$ cannot be bounded on $B^n$. $\square$

\smallskip

Now we can apply explicit calculations of $\o_d(Z)$ in Section \ref{entr.bounds} above to get explicit inequalities relating the geometry
of $Z$, the values of $f$ on this set, and the bounds on the derivatives of $f$. We shall restrict ourselves to the case of $Z$ being a
curve in the plane, as considered in Proposition \ref{curve}. Other situations presented in Proposition \ref{prop.o} can be treated in the
same way. Let $S$ be a connected smooth curve in $B^2$ of the length $\sigma$, and the injectivity radius $\e_0$. For
$d \leq {\sigma \over {24}}-1$ put $\kappa_d={1\over {2l}}(1-{{24}\over m}),$ in notations of Proposition \ref{curve}.

\bp\label{curve.smooth}
Let $f: B^2 \rightarrow {\mathbb R}$ be a $k$ times continuously differentiable function on $B^2$, with
$M_l(f)=\max_{B^2} \Vert d^l f \Vert, \ l=0,1,\ldots,k,$ and $S\subset B^2$ be a curve with the length $\sigma$, and with the
injectivity radius $\e_0$. Put $L=\max_{x\in S}|f(x)|$. Then for each $s \leq {\sigma \over {24}}-1$ we have

\be\label{curve.smooth.eq}
M_0(f)\leq ({8\over{\kappa_s}})^s(L+E^T_s(f))+E^T_s(f),
\ee
with $E^T_s(f)={1\over {(s+1)!}}M_{s+1}(f)$ and $\kappa_s={1\over {2l}}(1-{{24}\over m})>0.$
For each $s$ there are curves $S_s \subset B^2$ of the length at least $2s$ such that $M_0(f)$ cannot be bounded in terms of $L$ and
$M_l(f), \ l\geq s+1.$
\ep
\pr
The bound follows directly from Propositions \ref{smallZ}, \ref{curve}, and \ref{R.Const.Entr}. Now take as a curve $S_s$ a zero set
of a polynomial $y=T_s(x)$ inside $B^2$. Then for $f(x,y)=K(y-T_s(x))$ vanishing on $S_s$ \ $M_0(f)$ cannot be bounded through
$M_l(f), \ l\geq s+1$. $\square$

\medskip

Another way to extract more explicit answer from Proposition \ref{Taylor1} is to bound the norms $M_l(f)$ of the $l$-th order derivatives
of $f$, for $l=0,1,\ldots,k,$ by their maximal value $M = M(f)$, to substitute $M$ instead of $M_l(f)$ into the inequality
\ref{Smooth.Remez.ineq1}, and to explicitly minimize the resulting expression in $d$.

We shall fix the smoothness $k$ and consider sets $Z\subset B^n$ for which $\o(Z)=\o_{k-1}(Z)>0$. In particular, let $Z\subset B^n$ be
a measurable set with $m_n(Z)>0$. Then $\o_d(Z)\geq m_n(Z)$ for each $d$. Sets $Z$ in the specific classes, discussed in Section
\ref{entr.bounds} above, provide additional examples. Since $\o_0(Z)\geq \o_1(Z) \geq ... \geq \o_{k-1}(Z)$, by Proposition \ref{R.Const.Entr}
for each $d=0,\ldots,k-1$ we have $R_d(Z)\leq ({4n\over {\o(Z)}})^d$. Let us denote ${4n\over {\o(Z)}}\geq 4n$ by $q=q(Z)$.

\smallskip

The following theorem provides one of possible forms of an explicit inequality, generalizing the Remez one to smooth functions:

\bt\label{gen.bd}

Let $f: B^n \rightarrow {\mathbb R}$ be a $k$ times continuously differentiable function on $B^n$, with
$M_l(f)=\max_{B^n} \Vert d^l f \Vert \leq M=M(f), \ l=0,1,\ldots,k,$ and let a subset $Z\subset B^n$ with $\o_{k-1}(Z)>0$
be given. Put $L=\max_{x\in Z}|f(x)|$, $q=q(Z)\geq 4n$. Then

\be\label{Smooth.Remez.ineq2}
M_0(f)=\max_{x\in B^n} \vert f(x) \vert \leq 2q^{d_0}L+{1\over {(d_0+1)!}}M,
\ee
where $d_0=d_0(M,L),$ satisfying $1\leq d_0\leq k-1$, is defined as follows: $d_0=0$ if $L>M$, \ $d_0=k-1$ if
$L \leq {1\over {k!}}M$, and for ${1\over {k!}}M\leq L\leq M$ the degree $d_0$ is defined by 
${1\over {(d_0+1)!}}M\leq L \leq {1\over {d_0!}}M.$

\medskip

In particular, for $L>M$ the inequality takes the form

\be\label{dzero}
M_0(f)\leq L+2M,
\ee
while for $L \leq {1\over {k!}}M$ we get

\be\label{dzerok}
M_0(f)\leq 2q^{k-1} L+{1\over {k!}}M.
\ee
\et
\pr
As above, $R_d(Z)\leq ({4n\over {m_n(Z)}})^d=q^d.$ By Theorem \ref{Smooth.Remez.ineq} we have
$$
\max_{x\in B^n} \vert f(x) \vert \leq \inf_{d=0,1,\ldots,k} [q^d(L+E^T_d(f))+E^T_d(f)]\leq
$$
$$
\leq q^d(L+{1\over {(d+1)!}}M)+{1\over {(d+1)!}}M.
$$
Now we guess the value of $d$ which approximately minimizes the expression in the right-hand side: let $d_0=d_0(M,L)$ be defined as follows:

\smallskip

\noindent $d_0=0$ if $L>M$, \ $d_0=k-1$ if $L \leq {1\over {k!}}M$, and for ${1\over {k!}}M\leq L\leq M$ the degree $d_0$ is uniquely
defined by the condition
$$
{1\over {(d_0+1)!}}M\leq L \leq {1\over {d_0!}}M.
$$
In each case we have $1\leq d_0\leq k-1$. Substituting $d_0$ into the above expression we obtain for $L>M$ the inequality
$M_0(f)=\max_{x\in B^n} \vert f(x) \vert \leq L+2M$, while for $L\leq M$ we get $M_0(f)\leq 2q^{d_0} L+{1\over {(d_0+1)!}}M.$ In the case
$L \leq {1\over {k!}}M$ we get $d_0=k-1$, and the inequality takes the form $M_0(f)\leq 2q^{k-1} L+{1\over {k!}}M.$ $\square$

\medskip

\noindent{\bf Remark} In the case $L>M$ in Theorem \ref{gen.bd} we have $d_0=0$ and the resulting inequality (\ref{dzero}) is rather
straightforward. Indeed, we take one point $x_0\in Z$. By the assumptions, $|f(x_0)|\leq L$, while $||df||\leq M$ on $B^n$. For each
$x\in B^n$ we have $||x-x_0||\leq 2$. Hence $|f(x)|\leq L+2M.$ However, for smaller $L$, i.e. for larger $d_0$ the result apparently cannot
be obtained by a similar direct calculation. Compare a discussion in the next section.

\section{Whitney Extension of Smooth Functions}\label{Whitney}
\setcounter{equation}{0}

There is a classical problem of Whitney (see \cite{Fef} and references therein) concerning extension of $C^k$-smooth functions from closed 
sets. Recently a major progress have been achieved in this problem. The following ``Finiteness Principle'' has been obtained, in its 
general form, by Ch. Fefferman in 2003: for a finite set $Z\subset B^n$ and for any real function $f$ on $Z$ denote by $||f||_{Z,k}$ the 
minimal $C^k$-norm of the $C^k$-extensions of $f$ to $B^n$.

\smallskip

\noindent{\it There are constants $N$ and $C$ depending
on $n$ and $k$ only, such that for any finite set $Z\subset B^n$ and for any real function $f$ on $Z$ we have
$||f||_{Z,k}\leq C \max_{\tilde Z}||f||_{\tilde Z,k}$, with $\tilde Z$ consisting of at most $N$ elements.}

\smallskip

The original proof of this result, as well as its further developments in \cite{Fef} and other publications, provide rich
connections between the geometry of $Z$ and the behavior of the $C^k$-extensions of $F$. Effective algorithms for the extension have
been also investigated in \cite{Fef}. Still, the problem of an explicit connecting the geometry of $Z$, the behavior of $f$ on $Z$, and
the analytic properties of the $C^k$-extensions of $f$ to $B^n$ for $n\geq 2$ remains widely open. In one variable divided finite
differences provide a complete answer (Whitney). The following result illustrate the role of the Remez constant $R_d(Z)$ in the
extension problem.

\bt\label{Remez.Whitney}
For a finite set $Z\subset B^n$ and for any $x\in B^n\setminus Z$ let $Z_x=Z\cup \{x\}.$ Let $f_{Z,x}$ be zero on $Z$ and $1$ at $x$ and
let $\tilde f_{Z,x}$ be a $C^k$-extensions of $f_{Z,x}$ to $B^n$. Then for each $d=0,\ldots,k-1$ we have
$M_{d+1}(\tilde f_{Z,x})\geq {{(d+1)!}\over {R_d(Z)+1}}.$
\et
\pr
By Proposition \ref{Taylor1} we have for the extension $\tilde f_{Z,x}$
$$
M_0(\tilde f_{Z,x}) \leq \min_{d=0,1,\ldots,k-1} [R_d(Z)(L+E^T_d(f))+E^T_d(f)],
$$
where $E^T_d(\tilde f_{Z,x})={1\over {(d+1)!}}M_{d+1}(\tilde f_{Z,x})$ is the Taylor remainder term of $f$ of degree $d$ on the unit
ball $B^n$. In our case $M_0(\tilde f_{Z,x})\geq 1$ while $L=0$. So we obtain 
$1\leq \min_{d=0,1,\ldots,k-1} (R_d(Z)+1){1\over {(d+1)!}}M_{d+1}(\tilde f_{Z,x})$. We conclude that for each $d=0,\ldots,k-1$ we have 
$M_{d+1}(\tilde f_{Z,x})\geq {{(d+1)!}\over {R_d(Z)+1}}.$ $\square$

\smallskip

The results of Section \ref{Smooth.f} above can be translated into more results on extension from finite set, similar to that of Theorem
\ref{Remez.Whitney}. More important, Remez inequality for polynomials can be significantly improved, taking into account, in particular, a
specific position of $x$ with respect to $Z$. We plan to present these results separately.

\bibliographystyle{amsplain}

\end{document}